%% file: paper-ver-revised.tex
\documentclass[letterpaper, 10 pt, conference]{ieeeconf}  % Comment this line out if you need a4paper

\IEEEoverridecommandlockouts % This command is only needed if  you want to use the \thanks command

\usepackage{afterpage}
\usepackage{epsfig} % for postscript graphics files
\usepackage{times} % assumes new font selection scheme installed
\usepackage{amsmath} % assumes amsmath package installed
\usepackage{amssymb}  % assumes amsmath package installed
\usepackage{amsfonts}
\usepackage{mathptmx} % assumes new font selection scheme installed
\usepackage{tabularx}
\usepackage{epstopdf}
\usepackage{subfigure}
\usepackage{float}
\usepackage{epic,color}
\usepackage{mathrsfs}
\usepackage{dsfont,MnSymbol}
\usepackage{algorithm}
\usepackage{algorithmic}
\usepackage{multirow} % Multirow for tables
\usepackage[T1]{fontenc}

\newcounter{rmnum}

\newcounter{anum}

\def\IEEEQEDclosed{\mbox{\rule[0pt]{1.3ex}{1.3ex}}}

\def\qed{\ifmmode\IEEEQEDclosed\else{\unskip\nobreak\hfil
\penalty50\hskip1em\null\nobreak\hfil\IEEEQEDclosed
\parfillskip=0pt\finalhyphendemerits=0\endgraf}\fi}
\input{_definitions}

\def\eqms{\ \stackrel{L^2}{=}  \ }

\def\k{{\sf K}}
\def\barL{\mathcal{L}}
\def\barI{{\bar{I}}}
\def\err{{\cal E}}
\def\clU{\mathcal{U}}

\def\clF{{\cal F}}
\def\CrossTerm{\mathcal{C}}
\def\clN{\mathcal{N}}
\def\hax{\hat{x}}

\def\spmprev#1{} 
\def\pgmprev#1{}
\newlength{\noteWidth}
\long\def\notes#1{\ifinner
             {\tiny #1}
             \else
              \marginpar{\parbox[t]{\noteWidth}{\raggedright\tiny #1}}
               \fi}

\long\def\notes#1{}

\title{\LARGE \bf An Approach to Duality in Nonlinear Filtering}

\author{Jin-Won Kim, Amirhossein Taghvaei, Prashant G. Mehta and Sean
  P. Meyn% <-this % stops a space
\thanks{Financial support from the 
ARO grant W911NF1810334 and the NSF CMMI award 1761622 is gratefully acknowledged. 
}% <-this % stops a space
\thanks{J-W. Kim, A.~Taghvaei and P.~G.~Mehta are with the Coordinated
  Science Laboratory and the Department of Mechanical Science and
  Engineering at the University of Illinois at Urbana-Champaign
  (UIUC); S.~P.~Meyn is with the Department of Electrical and Computer
  Engineering at the University of Florida at Gainesville; Corresponding email: mehtapg@illinois.edu.}
}

\begin{document}

\maketitle
\thispagestyle{empty}
\pagestyle{empty}

%%%%%%%%%%%%%%%%%%%%%%%%%%%%%%%%%%%%%%%%%%%%%%%%%%%%%%%%%%%%%%%%%%%%%%%%%%%%%%%%
\begin{abstract}
	
This paper revisits the question of duality between minimum variance
estimation and optimal control first described for the linear Gaussian
case in the celebrated paper of Kalman
and Bucy.  A duality result is established for nonlinear filtering,
mirroring closely the original Kalman-Bucy duality of control and
estimation for linear systems. The result for the finite state-space
continuous time Markov chain is presented.  It's solution is used to
derive the classical Wonham filter.

% Nearly 60 years ago, in the celebrated paper of Kalman and Bucy, it was
% established that optimal estimation for linear Gaussian systems is dual to a
% linear-quadratic optimal control problem.  
% In this paper, for the first time, a duality result 
% The form of the result suggests a natural generalization which is
% presented as a conjecture for the continuous state
% space case.

\end{abstract}
%%%%%%%%%%%%%%%%%%%%%%%%%%%%%%%%%%%%%%%%%%%%%%%%%%%%%%%%%%%

\section{Introduction}

In Kalman's celebrated paper with Bucy,
it is shown that the problem of optimal estimation is dual
to an optimal control problem~\cite{kalman1961}.  A
striking example of the dual relationship is that, with the time arrow reversed,
the dynamic Riccati equation (DRE) of the optimal control is
the same as the covariance update equation of the Kalman filter.  
The relationship is useful, e.g., to derive results on
asymptotic stability of the linear filter based on asymptotic
properties of the solution of the DRE~\cite{ocone1996}.

A nonlinear extension of the minimum variance
estimator has been considered to be a harder problem.  In literature, 
it has been noted that: i) the dual relationship between the DRE of the
LQ optimal control and the covariance update equation of the Kalman
filter is {\em not} consistent with the interpretation of the negative
log-posterior as a value function; and ii) some of the linear
algebraic operations, e.g., the use of matrix transpose to define the
dual system, are not applicable to nonlinear
systems~\cite{todorov2006optimal,todorov2008general}.  
For these reasons, the original duality of Kalman-Bucy is seen as an LQG artifact
that does not generalize~\cite{todorov2006optimal}.  
 
In this paper, a nonlinear extension of the minimum variance estimation is
presented for the special case of a Markov process in
continuous time, on a finite state-space.  The dual system is a
backward ordinary differential equation.  An optimal control
objective is formulated whose solution yields the minimum variance
estimator.  Using the elementary method of change of control, the
formula for the optimal control is obtained and used to derive the
classical Wonham filter.

% \textit{The present   paper has a single contribution:} generalization
% of the original Kalman-Bucy duality theory to nonlinear filtering.
% It is an exact extension, in the sense that the dual optimal control
% problem has the same minimum variance structure for linear and nonlinear filtering
% problems.     The main result is developed for a Markov process in
% continuous time, on a finite state-space.  It is
% applied to derive the classical Wonham filter.  
% Extension to the continuous
% state-space is also examined: an optimal control problem is
% formulated and conjectured to be the dual in this setting.

\spm{Not sure we need or should emphasize this: The new duality theorem in this paper leads to an
elementary derivation of the classical nonlinear filter of Wonham.
}

%\smallskip

The outline of the paper is as follows: classical
duality is reviewed in \Sec{sec:prelim}, and the new dual optimal
control problem for the finite case is described in \Sec{sec:duality}.
Its solution leading to the Wonham filter is presented in \Sec{sec:main}.
\nobreak

%\newpage

\section{Background on classical duality}
\label{sec:prelim}

\newP{Linear Gaussian filtering model}   
Specified by
the linear
stochastic differential equation (SDE): 
%\begin{subequations}
	\begin{flalign}
&\text{Signal}\quad\quad\;\;&& \ud X_t =  A^\top X_t \ud t + \ud B_t&\nonumber\\%\label{eq:LG_dyn}\\
&\text{Observation}\;&& \ud Z_t = H^\top X_t\ud t + \ud W_t&\nonumber%\label{eq:LG_obs}
	\end{flalign}
%\end{subequations}
where $X_t \in \Re^d$ is the  state at time $t$, $Z_t \in \Re^m$ is the
observation, 
 $A$, $H$ are matrices of appropriate
dimension,
 and $B$, $W$ are mutually independent Wiener
processes (w.p.) taking values in $\Re^d$ and $\Re^m$,
respectively. The covariance matrices associated with $B$
and $W$ are denoted by ${Q}$ and ${R}$, respectively.  
The initial condition $X_0$ is drawn from a Gaussian
distribution $\clN(\hax_0,\Sigma_0)$, independent of $B$ or
$W$. It is assumed that the  noise covariance matrix is non-singular, ${R}\succ 0$.

% The filtering
% problem is to compute the posterior distribution $\P(X_t\mid \clZ_t)$
% where $\clZ_t:=\sigma(Z_\tau;\tau\in[0,t])$ denotes the filtration (the
% time-history of observations).  The solution to this problem is given
% by the well known equations of the Kalman filter. The duality arguments are briefly summarized in the following; cf.,~\cite[Ch. 7]{astrom1970}.

\newP{Minimum-variance estimator} 
Consider the problem of constructing a minimum variance estimator for the random variable $f^\top X_T$, at some fixed time $T$,  where $f\in\Re^d$ is an arbitrary, known vector. 

Given the observations $\{ Z_t :  t\in[0,T]\}$, the following linear
structure for the optimal estimator is assumed:
\begin{equation*}
S_T = y_0^\top \hax_0 - \int_0^T u_t^\top \ud Z_t
\label{eq:KF_LP_est}
\end{equation*}
where $y_0\in\Re^d$ is constructed below,  and the input $u=\{u_t :  t\in[0,T]\}$ is chosen to solve the optimization problem,
\[
\min_{u} \;\; \E (|S_T - f^\top X_T|^2)
\]  
The solution $S^*_T$ coincides with  the minimum-variance
estimator of $ f^\top X_T$.  

This stochastic optimization problem  is converted to a deterministic
optimal control problem via duality.

\newP{Dual optimal control problem} 
% The quadratic cost function is defined with respect to the two covariance matrices: 
% \begin{equation}\label{eq:Lagrangian_finite}
% \ell(y,u) := \half u^\top R u  + \half y^\top Q y
% \end{equation}
%The optimal control problem is:
\[
\begin{aligned}
\mathop{\text{Minimize}}_{u}   \quad & J(u) = \half\ y_0^\top
\Sigma_0 y_0 + \int_0^{T} \half u_t^\top R u_t  + \half y_t^\top Q y_t\ud t \\
\text{Subject to}  \quad & \frac{\ud y_t}{\ud t}
= -A y_t -H u_t,\quad y_T = f %\quad \lim_{t \to \infty}\pr_t(x) = \post (x)
\end{aligned}
\]
%\begin{flalign*}
%&\mathop{\text{Minimize}}_{u}   &&\quad J(u) = \frac{1}{2}\ y_0^\top
%\Sigma_0 y_0 + \int_0^{T} \ell(y_t,u_t) \ud t &\\
%&\text{Subject to}&&\quad \frac{\ud y_t}{\ud t}
% = -A y_t -H u_t,\quad y_T = f &%\quad \lim_{t \to \infty}\pr_t(x) = \post (x)
%\end{flalign*}
The process $\{y_t :  t\in[0,T]\}$ is referred to as the dual process.   The solution of
the optimal control problem yields the optimal control input, 
along with the vector $y_0$ that determines the minimum-variance
estimator  $S^*_T$. 

The Kalman filter is obtained by expressing $\{S^*_t(f) : t\ge 0,\ f\in\Re^d\}$ as the solution to a linear SDE~\cite[Ch.~7]{astrom1970}.      

\section{Duality for Nonlinear Filtering:\\ The Finite State space Case}\label{sec:duality}

\newP{Nonlinear filtering model} The finite state-space filtering problem is considered, in which the  state-space is the  canonical basis  $\mathbb{S} = \{e_1,e_2,\hdots,e_d\}$ in $\Re^d$.  

The Markovian state process $X=\{X_t :  t\in[0,T]\}$ evolves in continuous time,   taking values in $\mathbb{S}$.   This and the observation process $Z=\{Z_t :  t\in[0,T]\}$ are modeled by the SDE,
\begin{subequations}\label{eq:dyn-obs}
	\begin{flalign}
&\text{Signal}\quad\quad&&\;\; \ud X_t =  A^\top X_t \ud t + \ud B_t&\label{eq:dyn}
\\
&\text{Observation}  &&\;\; \ud Z_t = H^\top X_t\ud t + \ud W_t&\label{eq:obs}
	\end{flalign}
\end{subequations}
where $A\in\Re^{d\times d}$ is the \textit{rate matrix},   
 $H\in\Re^{d\times m}$,  
 $W$ is an $m$-dimensional w.p.\ with covariance $R\succ
0$. $B=\{B_t:t\in[0,T]\}$ is defined by
$$
B_t = X_t - \int_0^t A^\top X_\tau \ud \tau
$$
and it is a martingale since $A$ is the generator of the Markov process.
The initial distribution for  $X_0$ is  denoted
$\pi_0\in {\cal P}(\mathbb{S})$ where ${\cal P}(\mathbb{S})$ denotes the probability simplex in $\Re^d$. It is assumed that $X$, $W$ are mutually independent.  
  
%    \begin{remark}
  The linear observation model is chosen without loss of generality:   for any function $h \;\colon \; \mathbb{S}\to \Re$ we have
  $h(x) = H^\top x$ for $x\in \mathbb{S}$,  with $H_{i}=h(e_i)$.
%	The formulation is similar to the linear diffusion. However, it is because $X_t$ can have value on the canonical bases. In fact, the $i^{th}$ row of $H$ represents $h(e_i)$ where $h(\cdot)$ is an arbitrary function on $\mathbb{S}$.
%\end{remark}

Two  filtrations are required in this work:   $\clF = \{\clF_t : t\ge 0\}$ and
$\clZ = \{\clZ_t : t\ge 0\}$ where
\[
\clF_t:=\sigma(X_\tau,W_\tau  : 0\le \tau \le t)
\,,
\quad \clZ_t = \sigma(Z_\tau: 0\le \tau \le t)
\]
Let 
 $C_\clZ^p$ denote the family of $\Re^p$-valued,  continuous,
 and $\clZ$-adapted 
functions of time (the superscript ``$p$'' is omitted in the special case $p=1$).

The filtering problem is to compute the posterior distribution
$\P(X_t \in \varble \mid \clZ_t)$~\cite{bensoussan2018estimation}.   The solution is derived here through duality,  very much like in the classical linear setting.

\newP{The dual system}  A backward ordinary differential equation (ODE) on $\Re^d$,  
\begin{equation}\label{eq:dyn_y}
\frac{\ud Y_t}{\ud t} = -AY_t - HU_t,\quad Y_T = f
\end{equation}
whose  solution is
\[
Y_t = e^{A(T-t)} f + \int_t^T e^{A(\tau-t)}H U_\tau \ud \tau\,, \quad 0\le t\le T  
\]
% The matrices $\{ e^{At} : t \in [0,T]\}$ correspond to the Markov semi-group associated with the  rate matrix $A$.  

An optimal control problem is posed for the dual
system~\eqref{eq:dyn_y} whose solution yields the nonlinear filter.
This requires some restrictions on the class of control inputs.
The set of \textit{admissible control inputs} is defined as follows:
\begin{equation}\label{eq:admissible_control_defn}
\clU:=\left\{ U_t = \k_t^\top Y_t + V_t: \k \in
  C_\clZ^{d\times m}, 
  \; V \in
  C_\clZ^m , \; t\in[0,T]\right\}
\end{equation} 
We denote $U=\{U_t:t\in[0,T]\}$, $\k=\{\k_t:t\in[0,T]\}$ and
$V=\{V_t:t\in[0,T]\}$.  By construction, $\k$ and $V$ and
$\clZ$-adapted processes but $U$ may not be $\clZ$-adapted because of
the backward nature of the ODE~\eqref{eq:dyn_y}.   

% Note that $U_t$ is a affine feedback policy with $\clZ_t$ adapted coefficients. However, since the dynamics of $Y_t$ is defined backward-in-time, $U_t$, as well as $Y_t$, is not adapted to any forward filtration. 
The following proposition provides explicit representations for the
solution of the backward ODE~\eqref{eq:dyn_y}. Its proof appears in
  Appendix~\ref{apdx:pf_thm_Y_0}.  
\medskip

\begin{proposition}
\label{thm:Y_0}
Consider the backward ODE~\eqref{eq:dyn_y} with control input $U_t =
\k_t^\top Y_t + V_t$  where $\{\k_t:t\in[0,T]\}$ and $\{V_t:t\in[0,T]\}$ are given
$\clZ$-adapted processes.  
Then there exist  $\clZ$-adapted processes $\{\Phi_t, \eta_t, \kappa_t, \gamma_t:t\in[0,T]\}$, and $Y_0\in\clZ_T$, such that  for each $t\in[0,T]$, 
\[
Y_t = \Phi_t Y_0 + \eta_t, \quad
U_t = \kappa_t^\top Y_0 + \gamma_t
\]
%where $\Phi_t, \eta_t, \kappa_t, \gamma_t \in \clZ_t$ and $Y_0\in\clZ_T$.   
\end{proposition}
% \medskip
% The proof appears in the Appendix~\ref{apdx:pf_thm_Y_0}.
This proposition is used to define stochastic integral being used throughout the paper which is illustrated in the  Appendix~\ref{apdx:ito}.

\medskip

\newP{Minimum-variance estimator}    The problem of interest is
precisely as in the linear Gaussian case:  given a fixed time $T>0$, and $f\in\Re^d$, the goal is to obtain a representation for the minimum variance estimator for the random variable $f^\top X_T$.  

Given observations $Z=\{ Z_t : 0\le t\le T \}$  defined according to the model~\eqref{eq:obs},
the following linear structure for the estimator will be justified:  
\begin{equation}
S_T = Y_0^\top \pi_0 - \int_0^T U_t^\top \ud Z_t
\label{eq:NL_est}
\end{equation}
The vector $ Y_0$ is obtained from the solution
to~\eqref{eq:dyn_y}.     

The optimal control input is chosen as the solution to
the optimization problem:
\begin{equation*}\label{eq:min_var}
\min_{U\in\clU} \;\; \E [|S_T - f^\top X_T|^2]
\end{equation*}
Justification for the form \eqref{eq:NL_est} is provided through the
formulation of the dual control problem.

% and the resulting estimator is referred to as the minimum-variance
% estimator.
%\spmprev{Super confusing to me:
%In the duality based derivation
%of the nonlinear filter, the optimization problem is converted into a
%\emph{stochastic} optimal control problem that will be introduced
%after the following remark:}

\begin{remark}\label{rem:stoch_int}
The stochastic integral $\int_0^T U_t^\top \ud Z_t$
in~\eqref{eq:NL_est} is defined as a forward integral.  Formally, for
a given admissible choice of $\clZ$-adapted processes $\k$ and $V$, upon using the
representation in \Prop{thm:Y_0},
\[
\int_0^T U_t^\top \ud Z_t = Y_0^\top \int_0^T \kappa_t \ud Z_t  + \int_0^T
\gamma_t^\top \ud Z_t  
\]
where $\{\kappa_t:t\in[0,T]\}$, $\{\gamma_t:t\in[0,T]\}$ are adapted processes and therefore the associated
integrals are well-defined as standard  It\^o-integrals.  A self-contained
background on interpreting stochastic integrals for the \textit{non-adapted} processes
considered in this paper appears in Appendix~\ref{apdx:ito}.   
% Since $U_t$ is not adapted to $\clZ_t$ or ${\cal
%   F}_t$, the stochastic integral~\eqref{eq:NL_est} cannot be
% interpreted as a ordinary It\^{o} integral. Thus stochastic forward
% integral is defined yet in the similar way to the original one in the
\end{remark}

\pgmprev{Wir müssen wissen — wir werden wissen!}
\spm{Very cool, my poet friend!}

\newP{Dual optimal control problem}  
%The dual optimal control problem is:
\begin{subequations}\label{eq:opt-cont-finite}
\begin{align}
&\mathop{\text{Min}}_{U\in\clU} \ \ J(U) = {\sf E} \;\Big( \half |Y_0^\top X_0
- Y_0^\top \pi_0|^2 +  \int_0^T \half
U_t^\top R U_t \ud t \nonumber\\
&\quad+\int_0^T \half Y_t^\top \ud \langle X,X^\top \rangle_t Y_t +  \err_tU_t^\top \ud W_t + \err_t Y_t^\top \ud B_t\Big)\label{eq:opt-cont-finite-a}\\
&\text{Subject to} \ \ \frac{\ud Y_t}{\ud t} = -A Y_t - H
U_t,\quad Y_T = f \label{eq:opt-cont-finite-b} %\quad \lim_{t \to \infty}\pr_t(x) = \post (x)
\end{align}
\end{subequations}
where $\langle X,X^\top \rangle$ denotes the quadratic variation of the Markov
process $X$, and the {\em error
  process} $\err = \{\err_t:t\in[0,T]\}$ is defined as follows:
\begin{equation}\label{eq:et_defn}
\err_t := Y_0^\top(X_0-\pi_0) +\int_0^t U_\tau^\top \ud W_\tau +\int_0^t Y_\tau^\top \ud B_\tau
\end{equation}
As in Remark~\ref{rem:stoch_int}, the four stochastic integrals
appearing above are defined also as forward integrals (see
Appendix~\ref{apdx:ito}).  

The relationship between the optimal control objective $J(\cdot)$ and the minimum variance objective~\eqref{eq:min_var} is illustrated in the following proposition. The proof appears in the Appendix~\ref{apdx:opt_control}.
%In Appendix~\ref{apdx:opt_control}, the following identity is
%established for any arbitrary choice of an admissible control input:

\medskip
\begin{proposition}\label{prop:justification-of-cost}
Consider the state-observation model \eqref{eq:dyn-obs}, the linear estimator \eqref{eq:NL_est} and the dual optimal control problem \eqref{eq:opt-cont-finite}. For any arbitrary choice of an admissible control input,
\begin{equation*}
J(U) = \half \E [|S_T - f^\top X_T|^2]
\end{equation*}
This provides a justification for the objective function~\eqref{eq:opt-cont-finite-a} and moreover shows that $J(U)\ge 0$ for any admissible control. 
\end{proposition}

\medskip
\begin{remark}
Consider a deterministic control input of the form $U_t = k_t^\top Y_t + v_t$
where $\{k_t\}$, $\{v_t\}$ are deterministic functions of time (in particular, they
do not depend upon the observations).  Such a control is
trivially admissible. In this case, $\{Y_t\}$ is a deterministic function
of time and the error process $\err$ is a $\clF$-martingale. Consequently, 
$$
\E\Big(\int_0^T \err_tU_t^\top \ud W_t + \err_t Y_t^\top \ud B_t\Big) = 0
$$
and the objective function in~\eqref{eq:opt-cont-finite-a} simplifies to
$$
J(U) = \half Y_0^\top \Sigma_0Y_0 + \int_0^T \half U_t^\top R U_t +
\half Y_t^\top \E( Q(X_t) )Y_t \ud t
$$
where $\Sigma_0 := \E((X_0-\pi_0)(X_0-\pi_0)^\top)$ and $Q(\cdot)$ is a
$\mathbb{S}\to \Re^{d\times d}$ map defined as follows:
\begin{equation*}
Q(e_i) := \sum_{j\neq i}A_{ij}(e_j-e_i)(e_j-e_i)^\top,\quad i = 1,\ldots,d
\end{equation*}
The resulting problem is a deterministic LQ problem whose optimal solution
$\{U_t^*:t\in[0,T]\}$ will (in general) yield a sub-optimal estimate $S_T^*$
using~\eqref{eq:NL_est}.  The general problem considered here is much
tougher because $\err$ is {\em not} a $\clF$-martingale:  Under arbitrary admissible controls, it is not even adapted to this filtration.  
%If $\k_t$ and $V_t$ were deterministic, these terms vanish upon expectation as $U_t$, $Y_t$ become deterministic and $\err_t$ is ${\cal F}_t$-adapted. Remaining terms forms an exact analogy with the terminal condition and Lagrangian from the linear Gaussian case:
%$$
%\E \big(\half |Y_0^\top X_0 - Y_0^\top \pi_0|^2\big) = \half y_0^\top \E((X_0-\pi_0)(X_0-\pi_0)^\top)y_0
%$$
%and the running cost is
%\begin{align*}
%\E\Big(\int_0^T U_t^\top R U_t \ud t &+ \int_0^T Y_t^\top \ud\langle X,X^\top \rangle_tY_t\Big)\\
%&= \int_0^T u_t^\top R u_t + y_t^\top \E(Q(X_t))y_t \ud t
%\end{align*}
%where $Q(x)$ is a $\mathbb{S}\to \Re^{d\times d}$ map
%\begin{equation*}
%Q(e_i) := \sum_{j\neq i}A_{ij}(e_j-e_i)(e_j-e_i)^\top,\quad i = 1,\ldots,d
%\end{equation*}
\end{remark}
\medskip

We  have now set the stage   to derive the nonlinear filter via the
solution to the dual optimal control problem.

\section{Derivation of the Nonlinear Filter}
\label{sec:main}

Recall that an admissible input has the form $U_t = \k_t^\top Y_t +
V_t$ where $t\in [0,T]$. The goal is to obtain a formula for the gain process $\k = \{\k_t:t\in[0,T]\}$ such that the best choice of $V = \{V_t:t\in[0,T]\}$ is zero.  

This choice of input class can be regarded as an instance of the method of ``change of
control'' because $V$ represents the new variable for
control~\cite[Ch. 3.1]{bensoussan2018estimation}.

If $V_t\equiv 0$ then $\bar{Y}=\{\bar{Y}_t:t\in[0,T]\}$ solves the backward ODE
\[
\frac{\ud \bar{Y}_t}{\ud t} = -A\bar{Y}_t - H
\k_t^\top \bar{Y}_t,\quad \bar{Y}_T = f 
\]
and the associated control is denoted $\bar{U}_t = \k_t^\top
\bar{Y}_t$ for $t\in [0,T]$. With an arbitrary $V$,  the solution is expressed 
\begin{align*}
Y_t = \bar{Y}_t + \tilde{Y}_t,\quad 
U_t = \bar{U}_t + \tilde{U}_t
\end{align*}
where  $\tilde{Y}=\{\tilde{Y}_t: t\in [0,T]\}$ also solves a backward ODE:
\begin{equation}\label{eq:tilde_y}
\frac{\ud \tilde{Y}_t}{\ud t} = -A\tilde{Y}_t - H
\k_t^\top \tilde{Y}_t - H V_t,\quad \tilde{Y}_T = 0
\end{equation}
with $
\tilde{U}_t = \k_t^\top \tilde{Y}_t + V_t$ for $t\in [0,T]$.  

The error term is analogously split as $
\err_t = \bar{\err}_t + \tilde{\err}_t
$,
with
\[
\begin{aligned}
\bar{\err}_t &= \bar{Y}_0(X_0-\pi_0)+\int_0^t\bar{U}_\tau^\top \ud W_\tau + \int_0^t \bar{Y}_\tau^\top\ud B_\tau
\\
\tilde{\err}_t &= \tilde{Y}_0(X_0-\pi_0)+\int_0^t\tilde{U}_\tau^\top \ud W_\tau + \int_0^t \tilde{Y}_\tau^\top\ud B_\tau
\end{aligned}
\]
%It is noted that   $\tilde{U}_t = \tilde{Y}_t =
%\tilde{\err}_t = 0$ for all $t \in[0,T]$ in the case  $V_t\equiv 0$.

The optimal gain is described in the following theorem.
\begin{theorem}\label{thm:opt-soln}
Consider the optimal control problem \eqref{eq:opt-cont-finite}. 
For any non-zero $V \in C_\clZ^m$,
$$
J(U) \geq J(\bar{U})
$$
where the optimal gain is defined as following:
\begin{subequations}\label{eq:thm1}
\begin{align}
\ud \bar{\pi}_t = A^\top \bar{\pi}_t \ud t - \k_t^\top (\ud Z_t - H^\top \bar{\pi}_t \ud t),\quad \bar{\pi}_0 = \pi_0 \label{eq:thm1-a}\\
\k_t = - \E\big((X_t-\bar{\pi}_t)(X_t-\bar{\pi}_t)^\top |\clZ_t\big)HR^{-1},\quad t\in[0,T] \label{eq:thm1-b}
\end{align}
\end{subequations}
%Moreover, $\k$ with such property is uniquely determined.
\end{theorem}
%\begin{remark}
%It then follows that $\bar{U} = \{\bar{U}_t:t\in[0,T]\}$ is optimal, with the
%resulting optimal cost given by $J(\bar{U})= \half \E (|\bar{S}_T - f^\top X_T|^2)$ where $\bar{S}_T$ is the optimal estimate obtained using~\eqref{eq:NL_est} with $U=\bar{U}$ and $Y_0=\bar{Y}_0$.
%\end{remark}

\medskip
\subsection{Proof of Thm.~\ref{thm:opt-soln}}

It is simple calculation to see that
\[
J(U) =J(\bar{U})  + J(\tilde{U}) + \E(\CrossTerm)
\]
where the cross-term $\CrossTerm$ is defined by
%\begin{subequations}
%\begin{align}
%\CrossTerm &= \tilde{Y}_0^\top(X_0-\pi_0)(X_0-\pi_0)^\top\bar{Y}_0 \label{eq:cross-a}
%\\
%&+ \int_0^T \tilde{U}_t^\top R \bar{U}_t \ud t + \tilde{Y}_t^\top \ud \langle X,X^\top \rangle_t \bar{Y}_t \label{eq:cross-b}\\
%&+ \int_0^T(\tilde{\err}_t\bar{U}_t^\top+\bar{\err}_t\tilde{U}_t^\top) \ud W_t + \int_0^T (\tilde{\err}_t\bar{Y}_t^\top + \bar{\err}_t\tilde{Y}_t^\top) \ud B_t \label{eq:cross-c}
%\end{align}
%\end{subequations}
\begin{align*}
\CrossTerm &= \underbrace{\tilde{Y}_0^\top(X_0-\pi_0)(X_0-\pi_0)^\top\bar{Y}_0}_{\text{term (i)}}
\\
&+ \underbrace{\int_0^T \tilde{U}_t^\top R \bar{U}_t \ud t + \tilde{Y}_t^\top \ud \langle X,X^\top \rangle_t \bar{Y}_t}_{\text{term (ii)}} \\
&+ \underbrace{\int_0^T(\tilde{\err}_t\bar{U}_t^\top+\bar{\err}_t\tilde{U}_t^\top) \ud W_t + \int_0^T (\tilde{\err}_t\bar{Y}_t^\top + \bar{\err}_t\tilde{Y}_t^\top) \ud B_t}_{\text{term (iii)}}
\end{align*}
The strategy now is to choose $\k$ such that $\E(\CrossTerm)=0$ for all possible choices of  
$\clZ$-adapted $V$.  

\newP{Term (i)}  A standard technique of optimal
control theory dictates that the terminal condition term be expressed as an
integral by introducing a dual variable.
Towards this goal, we
introduce a vector-valued stochastic process $\bar{\pi} = \{\bar{\pi}_t:t\in[0,T]\}$ with
$\bar{\pi}_0=\pi_0$ (the prior).  At this point of time, we only
require that $\bar{\pi}$ is a $\clZ$-adapted process.  The dynamics of this process will be defined later.    

Using the process $\bar{\pi}$,  together with the requirement \eqref{eq:tilde_y} that $\tilde{Y}_T = 0$, we obtain
$$ 
\tilde{Y}_0^\top(\pi_0-X_0)(\pi_0-X_0)^\top\bar{Y}_0 = -\int_0^T \ud \big(\tilde{Y}_t^\top(\bar{\pi}_t-X_t)(\bar{\pi}_t-X_t)^\top\bar{Y}_t\big) 
$$
The differential is evaluated by an application of the product formula:\footnote{See Appendix~\ref{apdx:ito} for a justification of the product
formula for the class of (non-adapted) stochastic processes arising in this paper.} 
\begin{align*}
\ud&\big(\tilde{Y}_t^\top(\bar{\pi}_t-X_t)(\bar{\pi}_t-X_t)^\top\bar{Y}_t\big)\\
%=&(-A\tilde{Y}_t+H\k_t^\top\tilde{Y}_t+HV_t)^\top(\bar{\pi}_t-X_t)(\bar{\pi}_t-X_t)^\top\bar{Y}_t\ud t\\
%&+\tilde{Y}_t^\top\big((\ud\bar{\pi}_t-A^\top X_t\ud t - \ud B_t)(\bar{\pi}_t-X_t)^\top\\
%&\quad\quad+(\bar{\pi}_t-X_t)(\ud\bar{\pi}_t-A^\top X_t\ud t - \ud B_t)^\top \\
%&\quad\quad+ (\ud\bar{\pi}_t-A^\top X_t\ud t - \ud B_t)\\
%&\quad\quad\;\cdot(\ud\bar{\pi}_t-A^\top X_t\ud t - \ud B_t)^\top\big)\bar{Y}_t\\
%&+\tilde{Y}_t^\top(\bar{\pi}_t-X_t)(\bar{\pi}_t-X_t)^\top(-A\bar{Y}_t+H\k_t^\top\bar{Y}_t)\ud t\\
=&\,
\tilde{Y}_t^\top \Big\{  
	\big(\ud\bar{\pi}_t - A^\top \bar{\pi}_t\ud t+\k_tH^\top (X_t-\bar{\pi}_t)\ud t-\ud B_t\big)(\bar{\pi}_t-X_t)^\top\\
&+(\bar{\pi}_t-X_t)\big(\ud\bar{\pi}_t - A^\top \bar{\pi}_t\ud t+\k_tH^\top (X_t-\bar{\pi}_t)\ud t-\ud B_t\big)^\top\\
&+\ud\langle(\bar{\pi}-X),(\bar{\pi}-X)^\top\rangle_t
%(\ud\bar{\pi}_t-A^\top X_t\ud t - \ud B_t)(\ud\bar{\pi}_t-A^\top X_t\ud t - \ud B_t)^\top
\Big\}
	\bar{Y}_t-V_t^\top H^\top (X_t-\bar{\pi}_t)(X_t-\bar{\pi}_t)^\top\bar{Y}_t\ud t
\end{align*}
where $\langle(\bar{\pi}-X),(\bar{\pi}-X)^\top \rangle$ denotes
the quadratic variation of the process $\bar{\pi}-X$.  It is noted
that each of the term in the integral is a quadratic either in
$\tilde{Y}_t$ and $\bar{Y}_t$ or in $V_t$ and $\bar{Y}_t$.  

\newP{Term (ii)} The second term is expressed as:
%easy?   is the easiest to deal with.  It is expressed as
\begin{align*}
\int_0^T& \tilde{U}_t^\top R \bar{U}_t \ud t + \tilde{Y}_t^\top \ud \langle X,X^\top \rangle_t \bar{Y}_t \\
  &= \int_0^T \Big( \;\tilde{Y}_t^\top
  \big(\k_tR\k_t^\top \ud t +  \ud \langle X,X^\top \rangle_t\big)\bar{Y}_t + V_t^\top R\k_t^\top\bar{Y}_t\ud t \Big)
\end{align*}

\newP{Term (iii)} It remains to tackle the two stochastic integrals
involving the error processes. We begin by recalling~\eqref{eq:et_defn}:
\jin{I think this $\err_t$ is fine.}
\begin{align*}
\err_t &= Y_0^\top (X_0 - \pi_0) + \int_0^t U_\tau^\top \ud W_\tau +\int_0^t Y_\tau^\top \ud B_\tau
\end{align*}
Proceeding as in term~(i), the    process $\bar{\pi}$ is again used to express 
the terminal condition term $Y_0^\top (\pi_0 - X_0)$ as an integral.  Once again, using the product rule
\begin{align*}
\ud \big({Y}_t^\top  (X_t-\bar{\pi}_t)\big) 
%= (-AY_t-H\k_t^\top Y_t-HV_t)^\top(X_t-\bar{\pi}_t)\ud t+{y}_t^\top\ud (X_t-\bar{\pi}_t)\\
%=& {y}_t^\top(-A^\top +K_tH^\top)(\bar{\pi}_t-X_t)\ud t + V_tH^\top(\bar{\pi}_t-X_t)\ud t\\
%&+{y}_t^\top(\ud\bar{\pi}_t-A^\top X_t\ud t - \ud B_t)\\
=-Y_t^\top\big(\ud\bar{\pi}_t & - A^\top \bar{\pi}_t\ud t+\k_tH^\top (X_t-\bar{\pi}_t)\ud t\big)\\
& +  Y_t^\top\ud B_t -V_t^\top H^\top (X_t-\bar{\pi}_t)\ud t
\end{align*}
Therefore,
\begin{align*}
\err_t =& Y_0^\top(X_0-\pi_0) + \int_0^t U_\tau^\top \ud W_\tau + \int_0^t Y_\tau^\top \ud B_\tau\\
%=& Y_t^\top(\bar{\pi}_t-X_t) \\
%&-\int_0^t Y_\tau^\top\big(\ud\bar{\pi}_\tau - A^\top \pi_\tau\ud \tau-\k_\tau H^\top (X_\tau-\bar{\pi}_\tau)\ud \tau-\ud B_\tau\big) \\
%&-\int_0^tV_\tau^\top H^\top (X-\bar{\pi}_\tau)\ud \tau + \int_0^tY_\tau^\top \k_\tau + V_\tau \ud W_\tau\\
%&-\int_0^t Y_\tau^\top \ud B_\tau\\
%=& Y_t^\top(\bar{\pi}_t-X_t)-\int_0^t Y_\tau^\top\big(\ud\bar{\pi}_\tau - A^\top \pi_\tau\ud \tau-\k_\tau H^\top (X_\tau-\bar{\pi}_\tau)\ud \tau\big) \\
%&-\int_0^tV_\tau^\top H^\top (X-\bar{\pi}_\tau)\ud \tau + \int_0^tY_\tau^\top \k_\tau + V_\tau^\top \ud W_\tau\\
=&Y_t^\top (X_t-\bar{\pi}_t)+\int_0^t V_\tau^\top(\ud W_\tau + H^\top( X_\tau-\bar{\pi}_\tau)\ud \tau) \\
&+\int_0^tY_\tau^\top\big(\ud\bar{\pi}_\tau - A^\top \pi_\tau\ud \tau+\k_\tau(\ud W_\tau +H^\top(X_\tau-\bar{\pi}_\tau)\ud \tau)\big)
\end{align*}

In order to reduce the  notational burden,   the following
differential notation is adopted for the  $\clZ$-adapted stochastic processes
$\barI=\{\barI_t:t\in[0,T]\}$ and $\barL=\{\barL_t:t\in[0,T]\}$:
\begin{align*}
\ud \barI_t&:=
%\ud W_t +H^\top(X_t-\bar{\pi}_t)\ud t =
\ud Z_t -  H^\top\bar{\pi}_t \ud t\\
\ud \barL_t&:= \ud \bar{\pi}_t - A^\top \bar{\pi}_t\ud t + \k_t \ud \barI_t
\end{align*}
%It is noted that both $Z_t$ and $\bar{\pi}_t$ are forward $\clZ$-adapted.
The notation is used to express the error succinctly as
%$$
%\err_t =Y_t^\top (\bar{\pi}_t-X_t) - \int_0^tY_\tau^\top\big(\ud\bar{\pi}_\tau - A^\top \bar{\pi}_\tau\ud \tau-\k_\tau\ud \barI_\tau\big)+\int_0^t V_\tau^\top \ud \barI_\tau
%$$
$$
\err_t =Y_t^\top (\bar{\pi}_t-X_t) - \int_0^tY_\tau^\top\ud \barL_\tau+\int_0^t V_\tau^\top \ud \barI_\tau
$$
In particular, upon splitting $\err_t = \bar{\err}_t+\tilde{\err}_t$, we have
\begin{align*}
\bar{\err}_t &=\bar{Y}_t^\top (X_t-\bar{\pi}_t) + \int_0^t\bar{Y}_\tau^\top\ud \barL_\tau\\
\tilde{\err}_t &=\tilde{Y}_t^\top (X_t-\bar{\pi}_t) + \int_0^t\tilde{Y}_\tau^\top\ud \barL_\tau+\int_0^t V_\tau^\top\ud \barI_\tau
\end{align*}
We thus obtain a useful expression for term (iii):
{\small
\begin{align*}
&\int_0^T\tilde{\err}_t\bar{U}_t^\top+\bar{\err}_t\tilde{U}_t^\top \ud W_t + \int_0^T (\tilde{\err}_t\bar{Y}_t^\top + \bar{\err}_t\tilde{Y}_t^\top) \ud B_t\\
&=\int_0^T\Big\{
\tilde{Y}_t^\top(X_t-\bar{\pi}_t)\bar{Y}_t^\top\k_t
+\bar{Y}_t^\top (X_t-\bar{\pi}_t)\tilde{Y}_t^\top\k_t
+\big(\int_0^t\tilde{Y}_\tau^\top\ud\barL_\tau\big)\bar{Y}_t^\top \k_t \\
&\quad\quad+\big(\int_0^t\bar{Y}_\tau^\top\ud \barL_\tau\big)\tilde{Y}_t^\top\k_t +\big(\int_0^t V_\tau^\top\ud \barI_\tau\big)\bar{Y}_t^\top\k_t + V_t^\top\bar{Y}_t^\top (X_t-\bar{\pi}_t)\\
&\quad\quad+V_t^\top \big(\int_0^t\bar{Y}_\tau^\top\ud \barL_\tau\big)\Big\} \ud W_t
\\
&\quad+\int_0^T\Big\{\tilde{Y}_t^\top (X_t-\bar{\pi}_t)\bar{Y}_t^\top + \bar{Y}_t^\top (X_t-\bar{\pi}_t)\tilde{Y}_t^\top+\big(\int_0^t\tilde{Y}_\tau^\top\ud \barL_\tau\big)\bar{Y}_t^\top\\
&\quad\quad+\big(\int_0^t\bar{Y}_\tau^\top\ud \barL_\tau\big)\tilde{Y}_t^\top+\big(\int_0^t V_\tau^\top\ud \barI_\tau\big)\bar{Y}_t^\top\Big\}\ud B_t
\end{align*}}

This concludes our program of expressing each of three terms in
$\CrossTerm$ as an integral with sub-terms containing $\bar{Y}_t,\tilde{Y}_t,V_t$.  Now,
every sub-term is a quadratic of one of the two types: 
\begin{enumerate}
	\item The type 1 quadratic sub-terms contain $\bar{Y}_t$ and $\tilde{Y}_t$.  An
	example of this type of quadratic is $\tilde{Y}_t^\top
	\k_t R \k_t^\top\bar{Y}_t$ in the term (ii).
	\item The type 2 quadratic sub-terms contain $\bar{Y}_t$ and $V_t$.  An
	example of this is $V_t\k_t^\top\bar{Y}_t$ in the term (ii).
\end{enumerate}
We express
$
\CrossTerm = \CrossTerm_1 + \CrossTerm_2
$,
where $\CrossTerm_1$ contains only the quadratic sub-terms of type 1 and $\CrossTerm_2$ contains only
the quadratic sub-terms of type 2.  Upon collecting terms, we obtain
\spm{To do:  something seems wrong with the big equation here.   The two terms in the integrand at the end seem identical.}
{\small
\begin{align*}
\CrossTerm_1&=\int_0^T \tilde{Y}_t^\top (\k_tR\k_t^\top \ud t +  \ud \langle X,X^\top \rangle_t)\bar{Y}_t- \tilde{Y}_t^\top \ud\big\langle(\bar{\pi}-X_t),(\bar{\pi}-X_t)^\top\big\rangle_t\bar{Y}_t\\
&+\int_0^T\Big(\tilde{Y}_t^\top(\bar{\pi}_t-X_t)\bar{Y}_t^\top + \bar{Y}_t^\top(\bar{\pi}_t-X_t)\tilde{Y}_t^\top\Big)\ud\barL_t \\
%&+\int_0^T\tilde{Y}_t^\top\big(\ud\barL_t(\bar{\pi}_t-X_t)^\top + (\bar{\pi}_t-X_t)\ud\barL_t^\top\big)\bar{Y}_t \\
&+\int_0^T \Big(\big(\int_0^t \bar{Y}_\tau^\top\ud\barL_\tau\big)\tilde{Y}_t^\top \k_t
+ \big(\int_0^t \tilde{Y}_\tau^\top\ud\barL_\tau\big)\bar{Y}_t^\top \k_t\Big)\ud W_t\\
&+\int_0^T \Big(\big(\int_0^t \bar{Y}_\tau^\top\ud\barL_\tau\big)\tilde{Y}_t^\top
+\big(\int_0^t \tilde{Y}_\tau^\top\ud\barL_\tau\big)\bar{Y}_t^\top
\Big)\ud B_t
\end{align*} }
and {\small
\begin{align*}
&\CrossTerm_2 = \int_0^T V_t^\top\big(R \k_t^\top +H^\top( X_t-\bar{\pi}_t)(X_t-\bar{\pi}_t)^\top\big)\bar{Y}_t \ud t \\
&+\int_0^T \Big\{V_t^\top(X_t-\bar{\pi}_t)^\top \bar{Y}_t
+ V_t^\top \big(\int_0^t \bar{Y}_\tau^\top\ud\barL_\tau\big) 
+ \big(\int_0^t V_\tau^\top \ud \barI_\tau\big) \bar{Y}_t^\top \k_t\Big\} \ud W_t
\\
&\quad\quad +\int_0^T\big(\int_0^t V_\tau^\top \ud \barI_\tau\big)\bar{Y}_t^\top \ud B_t
\end{align*}}
In order to have $\E(\CrossTerm) = \E(\CrossTerm_1) + \E(\CrossTerm_2) = 0$ for
all possible choices of
$\bar{Y},\tilde{Y}$ and for all possible choices of $\bar{Y},V$, we follow the following 2-step procedure:
\begin{enumerate}
	\item In Step 1, we obtain an equation for $\bar{\pi}$ by setting 
	\[
	\E(\CrossTerm_1) = 0,\quad \text{a.s.}
	\]
	\item Given $\bar{\pi}$
	from Step 1,  we next derive a formula for the optimal gain
        $\k$ by imposing the requirement 
	\[
	\E(\CrossTerm_2) = 0,\quad \forall V \in C_\clZ^m
	\]  
\end{enumerate}
The 2-step procedure is   inspired by the analogous procedure in
classical LQ theory where the step 1 is used to derive the Ricatti equation and the
step 2 is used to derive the formula for the optimal feedback gain; cf.,~\cite[Ch. 7.3.1]{bensoussan2018estimation}.  

\newP{Step 1} By inspection, we find that upon setting
\begin{equation}
\ud \bar{\pi}_t = A^\top \bar{\pi}_t\ud t - \k_t \ud \barI_t \,,  \quad \bar{\pi}_0=\pi_0
 \label{eq:barpit}
\end{equation} 
which is as presented in the theorem statement~\eqref{eq:thm1-a}, we have
$\ud\barL_t \equiv 0$, and $\CrossTerm_1 $ reduces to
\begin{align*}
\CrossTerm_1=\int_0^T &\tilde{Y}_t^\top (\k_tR\k_t^\top \ud t+
\ud \langle X, X^\top \rangle_t)\bar{Y}_t
\\
-&\tilde{Y}_t^\top 
\ud\big\langle(\bar{\pi}-X),(\bar{\pi}-X)^\top\big\rangle_t\bar{Y}_t
\end{align*}
It is an easy calculation to compute the quadratic variation
\[
\ud \langle(\bar{\pi}-X),(\bar{\pi}-X)^\top\big\rangle_t = \k_tR\k_t^\top \ud t + \ud \langle X,X^\top \rangle_t
\]
and therefore, upon defining the dynamics of $\bar{\pi}$ according
to~\eqref{eq:barpit},  
\[
\CrossTerm_1  = 0  \quad \text{a.s.}
\]
%which is 0 that follows from the cross-variation terms computed due to Appendix~\ref{apdx:ito}.
This is true for {\em any} choice of $\clZ$-adapted gain process $\k$.  

\medskip

Among the  consequences are the following
pretty representations for the error processes: 
\begin{align}
\label{eq:error_rep_bar}
\bar{\err}_t &= \bar{Y}_t^\top (X_t-\bar{\pi}_t)+ \int_0^t
\bar{Y}_\tau^\top \ud\barL_\tau = \bar{Y}_t^\top (X_t-\bar{\pi}_t)
\end{align}
and similarly,
\begin{align*}
\tilde{\err}_t &= \tilde{Y}_t^\top (X_t-\bar{\pi}_t)+ \int_0^t V_\tau \ud\barI_\tau
\end{align*}
These expressions also hold for {\em any} $\clZ$-adapted $\k$.
% However, with the optimal choice of $\bar{k}_t$ (which remains to be
% determined), the process
% \[
% \bar{\pi}_t = \E(X\mid \clZ_t)
% \]
% The reasoning is as follows:
% \begin{enumerate}
% \item At time $t=0$, $\bar{\pi}_0=\E(X)$ because we chose
%   $\bar{\pi}_0=\pi_0$, the prior.
% \item At time $t=T$, $\bar{\pi}_T=\E(X\mid \clZ_T)$ because of how we set
%   the problem.
% \item At time $t\in(0,T)$, $\bar{\pi}_t=\E(X\mid \clZ_t)$ because of the
%   form of the admissible control.\footnote{I do not like this.  I wish
%     this conclusion was obtained as a mathematical formula!}
% \end{enumerate}
% Therefore, even though we do not yet know the form of the optimal
% $\k_t$, we know that $\bar{\pi}_t$ is a
% $\clZ_t$-martingale.\footnote{Alternatively, one can {\em assume} that
%   $\bar{\pi}_t$ is $\clZ_t$-martingale, derive the form of the optimal gain
%   $\k_t$ under this assumption and a posteriori justify the
%   assumption (see~Appendix~\ref{apdx:C}).} We derive
% the formula for $\k_t$ in the next step.  

\newP{Step 2} A formula for the
gain $\k=\{\k_t:t\in[0,T]\}$ is obtained by enforcing
$
\E [ \CrossTerm_2 ] = 0
$.

We first carry out some simplifications.
%In Appendix~\ref{apdx:calc},
It is straightforward calculation that, with $\bar{\pi}$ defined according
to~\eqref{eq:barpit}, the integrand of $\CrossTerm_2$ is a perfect
differential:  
\begin{equation}\label{eq:apdX_ref_1}
\CrossTerm_2 = \int_0^T \ud \big(\bar{\err}_t \int_0^t V_\tau^\top \ud
\barI_\tau\big) = \bar{\err}_T \int_0^T V_t^\top \ud
\barI_t
\end{equation}
The following orthogonality condition is thus obtained  upon using the representation~\eqref{eq:error_rep_bar} for $\bar{\err}_T$:
\spmprev{I don't think we can claim this is a projection theorem -- is my revision ok?}
\[
 f^\top \E \Big((\bar{\pi}_T - X_T) \int_0^T V_t^\top \ud \barI_t\Big)=\E (\CrossTerm_2) =0
\]
Since the function $f$ is arbitrary, we must have
\spm{To do (in October):   does the $\sigma$ algebra generated by $ \int_0^r V_t^\top \ud \barI_t$ ($V$ arbitrary, but $\clZ$ adapted) generate $\clZ_r$, for any $r$?  If so, we have Wonham}
\[
\E \Big((X_T-\bar{\pi}_T) \int_0^T V_t^\top \ud \barI_t\Big) = 0
\]  
To obtain the formula for $\k$, the expression inside the
expectation is written as an
integral---essentially by reversing the steps in finding the perfect
differential.  This
yields
\begin{equation}
\label{eq:apdX_ref_2}
\begin{aligned}
\E \Big(\int_0^T \big(\k_t &  R+(X_t-\bar{\pi}_t)(X_t-\bar{\pi}_t)^\top H\big)V_t\ud t\Big)  
\\
&-\;\; \E\Big(\int_0^T\big(\int_0^tV_\tau^\top \ud\barI_\tau\big) (\ud \bar{\pi}_t-A^\top X_t\ud t)\Big)= 0
\end{aligned}
\end{equation}
For the equation to hold for arbitrary choices of $V$ and
$\bar{I}$ (which is unrelated to the choice of
$V$), the two terms should both be zero:
\begin{align}
\E \Big(\int_0^T \big(\k_t R+(X_t-\bar{\pi}_t)(X_t-\bar{\pi}_t)^\top H\big)V_t\ud t\Big) &=0
\label{eq:step2_1}
\\
\E\Big(\int_0^T\big(\int_0^tV_\tau^\top \ud\barI_\tau\big) (\ud \bar{\pi}_t-A^\top X_t\ud t)\Big) &= 0
\label{eq:step2_2}
\end{align} 

The formula for the optimal $\k$ is obtained by
solving~\eqref{eq:step2_1}. 
% :
% \[
% E \Big(\int_0^T (\k_tR+(X_t-\bar{\pi}_t)(X_t-\bar{\pi}_t)^\top H)V_t \ud
% t \Big) = 0
% \]
Using the tower property of conditional expectation, because
$V_t$ and $\k_t$ are both $\clZ_t$-measurable, we have
\[
E \Big(\int_0^T (\k_tR+\E((X_t-\bar{\pi}_t)(X_t-\bar{\pi}_t)^\top H\mid \clZ_t) )V_t\ud
t \Big) = 0
\]
Since $V$ is an arbitrary $\clZ$-adapted function, $\k_t$ is uniquely determined on $L^2$ space: 
\[
\k_t=-\E((X_t-\bar{\pi}_t)(X_t-\bar{\pi}_t)^\top H\mid
\clZ_t)R^{-1},\quad t\in[0,T]
\]
This gives the formula for the optimal gain $\k$.  \qed

\begin{remark}
Using the optimal gain, the equation~\eqref{eq:barpit} for $\bar{\pi}$ becomes
\[
\ud \bar{\pi}_t = A^\top \bar{\pi}_t \ud t +\E [ (X_t-\bar{\pi}_t)(X_t-\bar{\pi}_t)^\top H\mid \clZ_t ]R^{-1}\ud
\barI_t,\quad \bar{\pi}_0 = \pi_0
\]
The equation is not closed because we do not know $\E(X_t\mid \clZ_t) =: \pi_t$.

One could consider closing the equation by
assuming a certainty equivalence principle that $\pi =
\bar{\pi}$.  In that case,
\[
\E((X_t-\bar{\pi}_t)(X_t-\bar{\pi}_t)^\top H\mid \clZ_t) = \text{diag}(\pi_t) (H - {\pi}_t^\top H)^\top
\]
where $\text{diag}(\pi_t)$ is a diagonal matrix whose diagonal entries are the elements of the vector $\pi_t$, and one obtains the equation
\[
\ud {\pi}_t = A^\top {\pi}_t \ud t +\text{diag}(\pi_t)(H - {\pi}_t^\top H)^\top  R^{-1}\ud
I_t,\quad {\pi}_0 = \pi_0
\]
where $\ud I_t = \ud Z_t - H^\top \pi_t \ud t$.  This is the equation for the
Wonham filter.
\end{remark}

\bibliographystyle{IEEEtran}
\bibliography{duality}

\appendix

\subsection{Proof of Proposition~\ref{thm:Y_0}}\label{apdx:pf_thm_Y_0}

For a given affine control law $U_t = \k_t^\top Y_t + V_t$, the
ODE~\eqref{eq:dyn_y} is a linear system with random coefficients:
\begin{equation}\label{eq:lin_ode}
\frac{\ud Y_t }{\ud t} = -(A+H\k_t^\top) Y_t - H V_t, \quad Y_T=f
\end{equation}
It admits a unique solution $Y:[0,T]\rightarrow \Re^d$. 
 Now,
because $\{ \k_t, V_t ; t\in[0,T]\}$ are $\clZ$-adapted and $Y_T=f$ is
deterministic, the solution $Y_0$  at time $t=0$ is a
$\clZ_T$-adapted random vector.     

For $t\geq\tau$, the state transition matrix $\Phi(t,\tau)$ is defined
as the solution to the matrix ODE
\begin{equation}\label{eq:transition_matrix}
\frac{\ud }{\ud t} \Phi(t,\tau) = -(A+H\k_t^\top) \Phi(t,\tau),\quad \Phi(\tau,\tau)=I
\end{equation}
A solution
of~\eqref{eq:lin_ode} is given by
\[
Y_t = \Phi(t\,;0) Y_0 - \int_0^t\Phi(t\,;\tau)HV_\tau \ud \tau =: \Phi_t Y_0 + \eta_t
\] 
%Note that $\Phi_t$ and $\eta_t$ are both $\clZ_t$ adapted and $Y_0$
%does not depend on time. In consequence, $U_t$ is expressed as
Similarly,
$\displaystyle
U_t = \k_t^\top Y_t + V_t = {\underbrace{(\Phi_t^\top \k_t)}_{\kappa_t}}^\top Y_0 + \underbrace{(\k_t^\top\eta_t + V_t)}_{\gamma_t}$.

% \begin{align*}
% U_t &= \k_t^\top Y_t + V_t\\
% &= \k_t^\top \Phi_tY_0 + \k_t^\top \eta_t + V_t\\
% &= (\Phi_t^\top \k_t)^\top Y_0 + (\k_t^\top\eta_t + V_t)\\
% &=: \kappa_t^\top Y_0 + \gamma_t
% \end{align*}\qed

\subsection{Stochastic integrals}\label{apdx:ito}

Recall the filtrations:
$\clF_t:=\sigma(X_0,B_\tau,W_\tau: 0\le \tau \le t)$ and $\clZ_t = \sigma(Z_\tau:
\tau\in[0,t])$,  $t\in [0,T]$. There are two types of stochastic processes:
\begin{enumerate}
	\item Adapted stochastic processes: $W,\:B,\:X\in \clF$ and 
	$Z,\:\bar{\pi},\:\bar{I},\:\k,\:V\in \clZ$.
	
	\item 
	Non-adapted stochastic processes: $Y,\:U,\:\err$,  and 
	their optimal and perturbed counterparts,
	$\bar{Y},\:\bar{U},\:\bar{\err}$ and
	$\tilde{Y},\:\tilde{U},\:\tilde{\err}$, respectively.  
\end{enumerate} 

Now, allowing only for admissible control inputs from $\clU$
(see~\eqref{eq:admissible_control_defn}), a generic stochastic process
considered in this paper is expressed as $\phi_t = F^\top \xi_t +
\alpha_t$,
 where $F\in \clZ_T$ and $\xi_t,\alpha_t\in \clF_t$ for each $t$
  (\Prop{thm:Y_0} and \Prop{prop:prop_error}).  

% \spmprev{yawn: For adapted processes, one may take $F=0$.    }

\medskip

\begin{definition}\label{defn:defn1}
	Consider two stochastic processes $\phi_t = F^\top \xi_t +
	\alpha_t$ and $\psi_t = G^\top \zeta_t + \beta_t$,  
	 where  $\xi_t,
	\alpha_t, \zeta_t , \beta_t\in {\cal F}_t$ are piecewise
        continuous functions of time $t$ with at most finitely many jumps and $F,\; G$ are bounded random vectors. Consider a partition
	$\Pi^N_{[0,t]} = \{0=t_0<t_1<\ldots<t_N = t\}$ with $\Delta:=\displaystyle\max_{i} (t_{i}-t_{i-1})$.  Then,
	\begin{align*}
	\int_0^t \phi_\tau \ud \psi_\tau & :=
	\lim_{\Delta \rightarrow 0} \sum_{i=1}^{N}\phi_{t_{i-1}}
	(\psi_{t_{i}}-\psi_{t_{i-1}})\\
	\langle \phi,\psi \rangle_t & :=
	 \lim_{\Delta
		\rightarrow 0}\sum_{i=1}^{N}
	(\phi_{t_{i}}-\phi_{t_{i-1}})  (\psi_{t_{i}}-\psi_{t_{i-1}})
	\end{align*}
provided the respective limits exist in $L^2$.
\end{definition}

\begin{proposition}\label{prop:prop2}
	Consider the two stochastic processes $\{\phi_t,\psi_t\}$ as defined
        in~Defn.~\ref{defn:defn1}. Then 
%the limit for the forward integral exists as
	\begin{align*}
	\int_0^t \phi_\tau \ud \psi_\tau  \;  \eqms &\; F^\top\Big(\int_0^t \xi_\tau \ud \zeta_\tau^\top\Big) G + F^\top\Big(\int_0^t\xi_\tau \ud \beta_\tau\Big)
\\
	& \quad \quad +\;\; G^\top\Big(\int_0^t \alpha_\tau \ud \zeta_\tau\Big) + \int_0^t \alpha_\tau \ud \beta_\tau\\
	% \end{align*}
	% and the limit exists for the cross variation as
	% \begin{align*}
	\langle \phi,\psi \rangle_t  \;\eqms &\; F^\top\langle \xi,\zeta^\top \rangle_tG + F^\top\langle \xi,\beta \rangle_t+\langle \alpha,\zeta^\top \rangle_tG + \langle \alpha,\beta \rangle_t
	\end{align*}
where the integrals on the right-hand side are standard 
It\^o-integrals.
%   \spmprev{the use of "forward Ito" might not clarify}
	% Where all integrals and quadratic variations above are defined using ordinary It\^{o} integral.
	% Therefore~Defn.~\ref{defn:defn1} is well defined.
	Moreover, the following It\^o product formula holds:
	\begin{equation*}
	\phi_t \psi_t - \phi_0 \psi_0 = \int_0^t \phi_\tau \ud \psi_\tau +
	\int_0^t \psi_\tau \ud \phi_\tau + \langle \phi,  \psi \rangle_t
	\end{equation*}
\end{proposition}
\medskip

\begin{proof}
The pre-limit is evaluated as	
\begin{align*}
	\sum_{i=1}^{N} &\phi_{t_{i-1}}(\psi_{t_{i}}-\psi_{t_{i-1}})
\\
%	=&\sum_{i=1}^{N} (F^\top\xi_{t_{i-1}}+\alpha_{t_{i-1}})(G^\top\zeta_{t_{i}}+\beta_{t_{i}} - G^\top\zeta_{t_{i-1}}-\beta_{t_{i-1}})\\
	%=&\sum_{i=1}^{N} F\xi_{t_{i-1}}(G\zeta_{t_{i}} - G\zeta_{t_{i-1}})
	%+ \sum_{i=1}^{N} F\xi_{t_{i-1}}(\beta_{t_{i}} - \beta_{t_{i-1}})\\
	%&+ \sum_{i=1}^{N} \alpha_{t_{i-1}}(G\zeta_{t_{i}}- G\zeta_{t_{i-1}})
	%+ \sum_{i=1}^{N} \alpha_{t_{i-1}}(\beta_{t_{i}} - \beta_{t_{i-1}})\\
	=&F^\top\Big(\sum_{i=1}^{N} \xi_{t_{i-1}}(\zeta_{t_{i}}^\top - \zeta_{t_{i-1}}^\top)\Big)G
	+ F^\top \sum_{i=1}^{N} \xi_{t_{i-1}}(\beta_{t_{i}} - \beta_{t_{i-1}})
\\
	&\quad+\; G^\top\sum_{i=1}^{N} \alpha_{t_{i-1}}(\zeta_{t_{i}}- \zeta_{t_{i-1}})
	+ \sum_{i=1}^{N} \alpha_{t_{i-1}}(\beta_{t_{i}} - \beta_{t_{i-1}})
	\end{align*}
The result is obtained upon letting $\Delta\to 0$.  For example,
%\spmprev{Ithink this is m.s., and not a.s.}
%	Taking limit $\Delta\to 0$ on each term leads the result. For the first term, the ordinary It\^{o} integral is defined by
$$
	\lim_{\Delta\to 0} \;\;\sum_{i=1}^{N} \xi_{t_{i-1}}(\zeta_{t_{i}}^\top - \zeta_{t_{i-1}}^\top)\;\;\eqms  \;\;\int_0^t\xi_\tau \ud \zeta_\tau^\top
$$
	and therefore, because $F,\;G$ are bounded, 
$$
	\lim_{\Delta\to 0} \;\;F^\top\Big(\sum_{i=1}^{N} \xi_{t_{i-1}}(\zeta_{t_{i}}^\top - \zeta_{t_{i-1}}^\top)\Big)G\;\;\eqms \;\; F^\top\Big(\int_0^t\xi_\tau \ud \zeta_\tau^\top\Big)G
$$
%where all limits are in  $L^2$. 
%	Similarly, the limits exist for the other terms.
	
The calculation for the cross variation is analogous.
%:  {\small 
%        \begin{align*} 
%	\sum_{i=1}^{N} &(\phi_{t_{i}}-\phi_{t_{i-1}})(\psi_{t_{i}}-\psi_{t_{i-1}})
%\\
%%	=&\sum_{i=1}^{N} (F^\top\xi_{t_{i}}+\alpha_{t_{i}} - F^\top\xi_{t_{i-1}}-\alpha_{t_{i-1}})(G^\top\zeta_{t_{i}}+\beta_{t_{i}} - G^\top\zeta_{t_{i-1}}-\beta_{t_{i-1}})\\
%	=&F^\top\Big(\sum_{i=1}^{N} (\xi_{t_{i}} - \xi_{t_{i-1}})(\zeta_{t_{i}}^\top - \zeta_{t_{i-1}}^\top)\Big)G
%	+ F^\top\sum_{i=1}^{N} (\xi_{t_{i}}-\xi_{t_{i-1}})(\beta_{t_{i}} - \beta_{t_{i-1}})
%\\
%	&+ \Big(\sum_{i=1}^{N} (\alpha_{t_{i}}-\alpha_{t_{i-1}})(\zeta_{t_{i}}^\top - \zeta_{t_{i-1}}^\top)\Big)G
%	+ \sum_{i=1}^{N}
%          (\alpha_{t_{i}}-\alpha_{t_{i-1}})(\beta_{t_{i}} -
%          \beta_{t_{i-1}})
%\\
%&\stackrel{\Delta\to 0}{\longrightarrow} F^\top\langle \xi,\zeta^\top \rangle_tG + F^\top\langle \xi,\beta \rangle_t+\langle \alpha,\zeta^\top \rangle_tG + \langle \alpha,\beta \rangle_t
%	\end{align*}}%
The product rule is proved by using the following identity
        (which holds for arbitrary stochastic processes): 
	\begin{align*}
	(\phi_{t_{i}}\psi_{t_{i}}-\phi_{t_{i-1}}\psi_{t_{i-1}})=&\phi_{t_{i-1}}(\psi_{t_{i}}-\psi_{t_{i-1}})+\psi_{t_{i-1}}(\phi_{t_{i}}-\phi_{t_{i-1}})\\
	&+(\phi_{t_{i}}-\phi_{t_{i-1}})(\psi_{t_{i}}-\psi_{t_{i-1}})
	\end{align*}
	Summing over $i$ and taking the limit as $\Delta\to 0$ yields
        the result.
\end{proof}

\begin{remark}
	The product rule is the {\em only} type of It\^o formula used in the various
	proofs in this paper.  This is because of the linear quadratic nature
	of the optimal control problem in finite-state-space settings.  The
	following differential notation is frequently used:
	\begin{equation}\label{eq:prod_rule}
	\ud (\phi_t \psi_t) = \phi_t \ud \psi_t + \psi_t \ud \phi_t + \ud
	\langle \phi,  \psi \rangle_t
	\end{equation}
\end{remark}

%\subsection{Relationship between the optimal
%control objective and the minimum variance
%objective}
\subsection{Proof of Proposition~\ref{prop:justification-of-cost}}
\label{apdx:opt_control}

The following identity is established
in this section for any admissible control:
% $\{U_t\}\in{\cal U}$
\[
J(U) = \half \E [|S_T - f^\top X_T|^2]
\]
The lefhand-side is the optimal
control objective as defined
in~\eqref{eq:opt-cont-finite-a}.  The righthand-side is the mean-squared
error. 
%  as defined
% in~\eqref{eq:min_var}
Recall that $S_T$ is the linear estimate as
defined by~\eqref{eq:NL_est}, $f\in\Re^d$ is deterministic, and $X_T$ is the 
hidden state at time $T$.

% At time $t=T$, the estimator is defined as (see~\eqref{eq:NL_est}): 
% \begin{equation}
% S_T = Y_0^\top \pi_0 - \int_0^T U_t^\top \ud Z_t
% \label{eq:NL_est_appdx}
% \end{equation}
% The minimum variance objective is 
% \[
% \min_{\clU} \;\; \E [|S_T - f^\top X_T|^2]
% \]  

The approach is to use the dual ODE~\eqref{eq:dyn_y} to express the
mean-squared error as an integral.  The product
formula~\eqref{eq:prod_rule} is used to obtain
\begin{align*}
\ud (Y_t^\top X_t) &= \ud Y_t^\top X_t + Y_t^\top \ud X_t + \ud\langle Y^\top \!\!\!,\;X\rangle_t
\\
&=(-Y_t^\top A^\top - U_t^\top H^\top) X_t \ud t + Y_t^\top (A^\top X_t\ud t + \ud B_t)
\\
&=-U_t^\top H^\top X_t \ud t + Y_t^\top \ud B_t
\end{align*}
which is shorthand for the integral equation
\[
Y_T^\top X_T  =  Y_0^\top X_0 + \int_0^T U_t^\top H^\top X_t \ud t + Y_t^\top \ud B_t
\]
With $Y_T=f$,  upon subtracting this equation
from~\eqref{eq:NL_est}, 
\[
f^\top X_T - S_T = (Y_0^\top X_0-Y_0^\top\pi_0) + \int_0^T U_t^\top \ud W_t + Y_t^\top \ud B_t
\]
With the definition of the error process $\err_t$ in~\eqref{eq:et_defn},
the  left-hand side is identified:  $f^\top X_T - S_T = \err_T$.  

%Since the objective function is to minimize $\E(\err_T^2)$, we use
%the 
The product formula~\eqref{eq:prod_rule} is then used to obtain
\[
\half \err_T^2 = \half \err_0^2 + \int_0^T \err_t \ud \err_t + \half \langle \err,\err \rangle_T
\]
The integral form~\eqref{eq:opt-cont-finite-a} of the objective
function follows from  evaluating each of the terms
as summarized in the following.

\begin{proposition} 
\label{prop:prop_error}
Consider the error process $\err=\{\err_t:t\in[0,T]\}$ defined
	in~\eqref{eq:et_defn}.  Suppose $U=\{U_t:t\in[0,T]\}$ is any admissible control.  Then  
	\begin{align*}
	\err_0^2 & = |Y_0^\top X_0-Y_0^\top\pi_0|^2
\\
	\int_0^T \err_t \ud \err_t & = \int_0^T \err_tU_t^\top \ud W_t +\int_0^T \err_t Y_t^\top \ud B_t 
\\
	\langle \err,\err\rangle_T & = \int_0^T U_t^2 \ud t + Y_t^\top \ud \langle X,X \rangle_t Y_t 
	\end{align*}
\end{proposition}

The proof is the direct application of \Proposition{thm:Y_0} and \Prop{prop:prop2}.

\end{document}

%% file: _definitions.tex
\def\qed{\hspace*{\fill}~\IEEEQED\par\endtrivlist\unskip}%\mbox{\rule[0pt]{1.3ex}{1.3ex}}}

\def\Re{\mathbb{R}}

\def\varble{\,\cdot\,}
\def\varble{\,\cdot\,}

\def\Proposition#1{Prop.~\ref{#1}}

\def\Sec#1{Sec.~\ref{#1}}

\def\notes#1{\marginpar{\tiny #1}\typeout{Notes!
Notes!
Notes!
}}
\renewcommand{\notes}[1]{\typeout{notes!}}
\def\FRAC#1#2#3{\genfrac{}{}{}{#1}{#2}{#3}}
\def\half{{\mathchoice{\FRAC{1}{1}{2}}%
{\FRAC{2}{1}{2}}%
{\FRAC{3}{1}{2}}%
{\FRAC{4}{1}{2}}}}

\def\Re{\field{R}}
\def\k{{\sf K}}

\def\Sec#1{Sec.~\ref{#1}}

\def\clZ{{\cal Z}}

\def\Sec#1{Sec~\ref{#1}}

\def\E{{\sf E}}

\def\spm#1{\notes{\textcolor{blue}{SPM: #1}}}

\def\jin#1{\notes{\textcolor{blue}{jin: #1}}}

\def\Sec#1{Sec.~\ref{#1}}

\def\P{{\sf P}}

\def\IEEEQEDclosed{\mbox{\rule[0pt]{1.3ex}{1.3ex}}}
\def\qed{\nobreak\hfill\IEEEQEDclosed}

\def\clZ{{\cal Z}}
\def\varble{\,\cdot\,}

\newtheorem{theorem}{Theorem}

\newtheorem{definition}{Definition}

\newtheorem{remark}{Remark}
\newtheorem{proposition}{Proposition}

\def\beq{\begin{eqnarray}} 
\def\bc{\begin{center}} 
\def\be{\begin{enumerate}}
\def\bi{\begin{itemize}} 
\def\bs{\begin{small}}
\def\bS{\begin{slide}}
\def\ec{\end{center}} 
\def\ee{\end{enumerate}}
\def\ei{\end{itemize}}
\def\es{\end{small}}
\def\eS{\end{slide}}
\def\eeq{\end{eqnarray}}

\newcommand{\newP}[1]{\medskip\noindent{\bf #1:}}

\newcommand{\ud}{\,\mathrm{d}}

\def\Re{\mathbb{R}}
\def\E{{\sf E}}

\def\varble{\,\cdot\,}

\def\Sec#1{Sec.~\ref{#1}}

\def\Prop#1{Prop.~\ref{#1}}

\def\clZ{{\cal Z}}

\renewcommand{\Re}{\mathbb{R}}

\def\FRAC#1#2#3{\genfrac{}{}{}{#1}{#2}{#3}}